\documentclass[a4paper,12pt]{article}
\usepackage{amssymb,latexsym,amsmath,amsthm}

\theoremstyle{definition}
\newtheorem*{opr}{Definition}
 \theoremstyle{plain}
\newtheorem*{thm}{Theorem}
\newtheorem{lem}{Lemma}[section]
\newtheorem{utv}{Proposition}[section]
\newtheorem*{sled}{Corollary}
\newtheorem*{main}{Main Lemma}
\theoremstyle{remark}

\newcommand{\bl}{\square}

\newcommand{\C}{\mathbb{C}}
\newcommand{\Pm}{\mathbb{P}}
\newcommand{\doc}{\emph{Proof. }}
\renewcommand{\leq}{\leqslant}
\renewcommand{\geq}{\geqslant}

\newcommand{\ph}{\varphi}
\newcommand{\itd}[2]{#1_1,\dotsc, #1_{#2}}

\newcommand{\La}{\Lambda}

\newcommand{\odo}[1]{=1,\dots, #1}

\begin{document}
 
  \title {On the dimensions of  commutative subalgebras and   subgroups of nilpotent algebras and Lie groups of class 2}

\date{}
\author{Maria V. Milentyeva\thanks{This research is partially supported by RFBR (no.~05-01-00895.)}} 
\maketitle
 
\begin{abstract} 
We obtain the functions that bound the dimensions of finite dimensional nilpotent associative or Lie algebras of class 2 over an algebraically closed field in terms of the dimensions of their commutative subalgebras. As a result, we also compute a similar function for complex nilpotent Lie groups of class 2.
\end{abstract}

\section{Introduction}
We shall consider the following functions, introduced by the author  (Milentyeva,~2006):
 \begin{opr}
For any integer $n$ we denote by $l^n_{K}(n)$ ($a^n_{K}(n)$ and $g^{n}_K(n)$) the greatest integer $h$ such that there exists a nilpotent Lie algebra (an associative algebra and a Lie  group respectively) of class 2 and of dimension $h$ over a field $K$ such that the dimensions of all commutative subalgebras (Lie subgroups) of this algebra (Lie group) are not greater than $n$. (In the case of Lie groups $K$ is the complex or real   field.)
\end{opr}

Throughout this paper we assume that all algebras are finite dimensional.

It was shown  in Theorem~8 of (Milentyeva,~2006) that these functions have quadratic growth. More precisely, if  $f$ is one of the functions  $l^n_K,\ a^n_K$ or $g_K^n$, then the following estimate holds:
\begin{equation} \label{A}
 \frac {n^2+4n-5}{8} \leq f(n)\leq \frac{n^2}4+n.
\end{equation}
 
Analogs of these functions for finite $p$-groups were studied by A.Yu.~Olshanskii~(Olshanskii,~1978). Similar functions for finitely generated torsion-free nilpotent groups   and for arbitrary associative  algebras,  Lie algebras and Lie groups were considered  by the author in (Milenteva, 2004) and (Milentyeva,~2006) respectively. It was proved that these functions also have quadratic growth. It had not been known before that these functions are non-linear. In particular, the paper (Milenteva,~2004) 
provides an answer to   Question~8.76 of the Kourovka Notebook~(Mazurov, Khukhro,~1992), raised in 1982 by John Wilson, on the relationship between the torsion-free rank of a finitely generated nilpotent group and the ranks of its abelian subgroups.
 
Here we compute the functions $l^n_K(n),\ a^n_K(n)$, and  $g_K^n(n)$  in the case of algebraically closed field~$K$.

\begin{thm}
If $K$ is an algebraically closed field, then the functions   defined above are given by the following formula:
\begin{equation} \label{E:m}
l^n_K(n)= a^n_K(n)= g_{\C}^n(n)=
\begin{cases}
2n-1, & 1\leq  n\leq 7;\\
\left[ \frac{n^2+4}8  \right]+n, &n\geq 8.
\end{cases}
\end{equation}
\end{thm}

The Theorem shows that  if  $K$ is an algebraically closed field, the lower bound  from (Milentyeva,~2006) is asymptotically exact. It is easy to modify the proof of the relation~(\ref{A}) to achieve the  lower bound we need. The proof of the upper bound equal to the lower one is more complicated. The crucial step of this proof is the Main Lemma, 
which is the converse of Lemma~11 of (Milentyeva, 2006), that provides the existence of alternating bilinear forms corresponding to "large" nilpotent Lie algebras of class~2 with "small" dimensions of abelian subalgebras.

The structure of the paper is as follows.
Section~2 is devoted to the   Main Lemma. 
In~\S2.1 we recall some preliminary definitions and results from algebraic geometry. In~\S2.2  the Main Lemma is proved.
In \S3.1 we show that the functions $l^n_K(n),\ a^n_K(n)$, and  $g_K^n(n)$ are equal. Thus we reduce the Theorem to the case of Lie algebras. Finely, in~\S3.2 we apply  the Main Lemma to compute $l^n_K(n) $.
  
\section{The Main Lemma}
 In this section, $K$ is a given algebraically closed field.

\begin{main}
Suppose that the positive integers $k, \ t$ and $n$ satisfy the inequalities    
\begin{align} 
 2n\geq t(k-1)+2k,&      \label{E:c}\\
  t\geq 2;& \label{E:t2}
\end{align}
and let $V$ be an $n$-dimensional vector space over  $K$. Then for any $t$-tuple  ${\varPhi=\{\ph_1, \dotsc, \ph_t \}}$ of alternating bilinear forms on $V$ there exists a $k$-dimensional subspace that is simultaneously isotropic for all of the forms $\ph_1, \dotsc, \ph_t$.   
\end{main}

\subsection{Preliminaries.}

 The proof of the Main Lemma is based on some ideas of algebraic geometry.
So, we  recall that $X \subset {\Pm}^n$ is a \emph{closed in Zariski topology subset} if it consists of all points at which a finite number of homogeneous  polynomials with coefficients in $K$ vanishes.
This topology induces Zariski topology on any subset of  the projective space ${\Pm}^n$. A closed subset of ${\Pm}^n$ is    a \emph{ projective variety}, and an open subset of a projective variety is   a \emph{quasiprojective variety}. A nonempty   set $X$ is called \emph{irreducible} if it cannot be written as the union of two proper closed subsets.

Further, let   $f:X\to\Pm$$^m$ be a map of a  quasiprojective variety $X\subset\Pm$$^n $ to a projective space $\Pm$$^m$. This map is \emph{regular} if for every point $x_0 \in X$ there exists a neighbourhood $U \ni x_0$ such that the map $f: U \to \Pm$$^m$ is given by an $(m+1)$-tuple $(F_0:\ldots:F_m)$ of homogeneous polynomials of the same degree in the homogeneous coordinates of $x \in {\Pm}^n$, and $F_i(x_0)\neq 0$ for at least one $i$.

We will  use the following properties of the \emph{dimension} of a quasiprojective variety.
\begin{enumerate}
\item [{\sffamily(i)}]
The dimension of $\Pm$$^n$ is equal to $n$.
\item [{\sffamily(ii)}]
If $X$ is an irreducible variety and $U\subset X$ is open, then 
  $\dim U=\dim X$. 
\item [{\sffamily(iii)}]
The dimension of a reducible variety is the maximum of the dimension of its irreducible components.
\end{enumerate}

For more details see (Shafarevich,~1994). We will need the following propositions which are also to be found in (Shafarevich,~1994).

\begin{utv} \label{T:a}
   The property that a subset $Y \subset X$ is closed in a quasiprojective variety $X$ is a local property. That is, if  ${X=\cup U_\alpha}$ with open sets $U_\alpha$, and   $Y\cap U_\alpha$ is closed in  $U_\alpha$ for each $U_\alpha$, then $Y$ is closed in $X$.
\end{utv}

\begin{utv} \label{T:b}
  The image of a projective variety under a regular map is closed.
\end{utv}

\begin{utv} \label{T:B}
If $X\subset Y$ then  $\dim X\leq \dim Y$. If $Y$ is irreducible and $X\subset Y$ is a closed subvariety with   ${\dim X= \dim Y}$ then $X=Y$.
\end{utv}

\begin{utv} \label{T:c}
  Let $f:X \to Y$ be a regular map between irreducible varieties. Suppose that $f$ is surjective: $f(X)=Y$, and that $\dim X=n, \ \dim Y =m$. Then $m \leq n$, and
 \begin{enumerate}
\item [{\sffamily(i)}]    $\dim F\geq n-m$ for any $y\in Y$ and for any component $F$ of the fibre $f^{-1}(y)$;
\item [{\sffamily(ii)}] there exists a nonempty open subset $U\subset Y$ such that $ \dim f^{-1}(y)=n-m$ for ${y\in U} $.
\end{enumerate}
\end{utv}

Let ${\mathbb P}^n ,\ {\mathbb P}^m$ be projective spaces having homogeneous coordinates ${(u_0: \dotsc: u_n)}$  and \linebreak 
 ${(v_0: \dotsc: v_m)} $ respectively. Then the set   ${\mathbb P}^n \times {\mathbb P}^m$ of pairs $(x,y)$ with $x\in {\Pm}^n$ and $y\in {\Pm}^m$ is naturally embedded as a closed set into the projective space ${\mathbb P}^N$ where $N= (n+1)(m+1) -1 $. Thus there is a topology on ${\mathbb P}^n \times {\mathbb P}^m$, induced by the Zariski topology on $\Pm$$^N$.
 
\begin{utv} \label{T:e}
  A subset $X \subset {\mathbb P}^n \times {\mathbb P}^m\subset \Pm$$^N$ is a closed algebraic subvariety if and only if it is given by a system of equations  
\[
 G_i(u_0: \dotsc: u_n;v_0: \dotsc: v_m)=0 \qquad (i \odo t)
\]
homogeneous in each set of variables $u_i$ and $v_j$.
\end{utv}

\begin{utv} \label{T:A}
Let  $f:X\to Y$ be a regular map between projective varieties, with ${f(X)=Y}$. Suppose that $Y$ is irreducible, and that all the fibres  $f^{-1}(y)$ for $y \in Y$  are irreducible and of the same dimension. Then $X$ is irreducible.
\end{utv}
\begin{utv} \label{T:u}
   Consider an  $n$-dimensional vector space  $V$ with a basis  $\{e_1, \dotsc ,e_n\}$. Let   $U$ be a  $k$-dimensional subspace of  $V$  with a basis $\{f_1, \dotsc ,f_k\}$. To $U$  we  assign the point  $P(U)$ of the projective space $\mathbb{P}$$(\Lambda^k V)$ by the rule 
$$P(U)={f_1\wedge \dotsb \wedge f_k}.$$
The point $P(U)$ has the following form in the basis $\{e_{i_1}\wedge \dotsb \wedge e_{i_k} \}_{i_1< \dotsb <i_k} $ of $\La^k V$:
$$P(U)= {\sum_{i_1< \dotsb <i_k}p_{i_1 \dotsc i_k} e_{i_1}\wedge \dotsb \wedge e_{i_k}}.$$
Then the homogeneous coordinates $p_{i_1 \dotsc i_k}$ of  $P(U)$ are called the Plucker coordinates of $U$,
$P(U)$ is uniquely determined by   $U$, and the following assertions hold.

\begin{enumerate} 
 \item[{\sffamily(i)}]    
The subset of all points  $p\in \mathbb{P}$$(\Lambda^k V)$ of the form $p=P(U)$   is  closed in  $\mathbb{P}$$(\Lambda^k V)$; this subset $G(k,n)$ (the  Grassmanian or Grassmann variety) is defined by the relations   
\begin{equation} \label{E:d}
      \sum_{r=1}^{k+1}(-1)^r p_{i_1 \dotsc i_{k-1} j_r}p_{j_1 \dotsc \widehat{j_r} \dotsc j_{k+1}}=0
\end{equation} 
for all sequences   $i_1\dotsc i_{k-1}$ and $j_1 \dotsc j_{k+1}$.

\item[{\sffamily(ii)}]
  $\dim G(k,n)=k(n-k)$.
\item [{\sffamily(iii)}]  $G(k,n)$ is irreducible {\normalfont(see  Section~2.2.7 of (Vinberg, Onishchik,~1995))}.

\item[{\sffamily(iv)}] Suppose, for example, that  $p_{1\dotsc k} \ne 0$. If  $p=(p_{i_1\dotsc i_r})=P(U)$, then  $U$ has a basis $\{ \itd f k \}$ such that   
\begin{equation} \label{E:e}
 f_i=e_i+\sum_{r>k} a_{ir}e_r \qquad \mbox{for }  i\odo k,
\end{equation}
where
\begin{equation} \label{E:f}
  a_{ir}=(-1)^{k-i} \frac{p_{1\dotsc \widehat{i}\dotsc kr}}{p_{1\dotsc k}}.
\end{equation}
(Here $\widehat i $ means that the index $i$ is discarded.)
\end{enumerate}

\end{utv}

\begin{utv} \label{T:v} 
{\normalfont (see Section~14.7 of (Fulton,~1984) or Section~1.5 of (Griffits, Harris,~1978)) }
Concider an $n$-dimensional vector space $V$. For any increasing sequence $$
0 \subsetneq V_1\subsetneq \dotsc \subsetneq V_k 
$$
of subspaces of $V$
 put
\begin{equation*}
W(\itd V k)=\{p=P(U)\in G(k,n) \mid \dim (U\cap V_i)\geq i  \mbox{ for } i\odo k \}.
\end{equation*}
Then all the subsets $W(\itd V k)$, called Schubert cells, are closed in $G(k,V)$, and  
\begin{equation} \label{E:Y}
\dim W(\itd V k)= \sum_{i=1}^k (a_i -i) 
\end{equation}
where $a_i=\dim V_i$.
\end{utv}

\begin{sled}
Let $V$ be a vector space of dimension $n$, and let ${s_0=\max \{ 0, 2k-n\}}$. Given a $k$-dimensional subspace ${U\subset V}$, put
\begin{equation} \label{E:Z}
 G_s=G_s(U,V)=\{P(U')\in G(k,n) \mid \dim U \cap U' \geq s  \} \ \  \mbox{ for } s=s_0, \dotsc, k.
\end{equation}
Then, for each $s$, the sets $G_s$ is closed in $G(k,n)$, and
 \begin{equation} \label{E:L}
\dim G_s =  (k-s)(n-k+s). 
\end{equation}
\end{sled}
\doc
It is easy to check that the sets $G_s$ are Schubert cells. Indeed, choose a basis $\{\itd e n\}$ of $V$ such that $\itd f k$ span $U$. Then    ${G_s=W(\itd V k)}$ where
\begin{alignat*}{2}
&     V_i= \langle \itd e {k-s+i}\rangle\ \ \    &&\mbox { for } i \odo s, \\
&V_i = \langle \itd e {n-k+i}\rangle   &&\mbox{ for } i=s+1,\dotsc, k.
\end{alignat*}
(The conditions $\dim(U'\cap V_i)\geq i$ for $i> s$ are trivial , and for $i\leq s$ they follows from the condition
 $\dim (U' \cap V_s)=\dim(U'\cap U) \geq s.$)  Using~(\ref{E:Y}), we get
$$
\dim G_s= \sum_{i=1}^s(k-s)+\sum_{i=s+1}^k(n-k)=s(k-s)+(k-s)(n-k)=(k-s)(n-k+s) .\ \bl
$$

\subsection{Proof of the Main Lemma.}
Let $V$ be a vector space of dimension $n$ over $K$. Consider the vector space
$$L=\{ (\ph_1,\dotsc,\ph_t) \mid \ph_i \mbox{ is an alternating bilinear form on }   {V}\} .$$ 
Let ${\mathbb P}(L)$ be the corresponding projective space. We have $\dim {\Pm}(L) = t\frac{n(n-1)}{2}-1$.
Put  ${M={\mathbb P}(\Lambda^k  {V}) \times {\mathbb P}(L)}$. 
We regard $M$ as a closed subset of ${\mathbb P}^m$ where $m=\frac{tn(n-1)}{2} (_{\,k}^{\,n}) -1$. Let $S$ be the subset of $M$ consisting of pairs   $(p,\ph)$ such that  
\begin{enumerate}
\item [{\sffamily(i)}] $p\in G(k,n)$, that is, $p=P(U)$ for some $k$-dimensional subset  $U \subset  {V}$;
\item [{\sffamily(ii)}]  $U$ is isotropic for all forms of the $t$-tuple $\ph=(\ph_1,\dotsc,\ph_t)$.
\end{enumerate}

\begin{lem} \label{T:p}
$S$ is a projective variety. 
\end{lem}
\doc
Let us show that $S$ is closed in ${\mathbb P}^m$. To prove this we concider the covering of the projective variety $M$ by the open sets ${A_{i_1 \dotsc i_k}=\{ (p,\ph)\in M | p_{i_1 \dotsc i_k} \ne 0\}}$ and verify that  $A_{i_1 \dotsc i_k} \cap S$ is closed in~$A_{i_1 \dotsc i_k}$ for any sequence  $i_1 \dotsc i_k $ with $ i_1 <\dotsb <i_k $. Then, by Proposition~\ref{T:a}, $S$ is closed in $M$. Hence, since $M$ is closed in ${\mathbb P}^m$, $S$ is also closed in ${\mathbb P}^m$.

The condition for $p\in G(k,n)$ is defined by the relations~(\ref{E:d}).

We can assume without loss of generality that $p_{1\dotsc k} \ne 0$. (The other cases are similar.) If $p=P(U)$, then $U$ has a basis of the form~(\ref{E:e}). The space $U$  is isotropic for all forms $\ph_l$ with 
  $ l\odo t$ if and only if
\[
 \ph_l(f_i,f_j)=0  \qquad  \mbox{for } l\odo t \mbox{ and } i,j\odo k.
\]             

We substitute for $f_i$ using~(\ref{E:e})   and then replace  $a_{ir}$  using~(\ref{E:f}). Now, multiplying both sides by $p_{1\dots k}^2$  and recalling that  ${\ph_l(e_i,e_j)=\ph^l_{ij}}$, we obtain a system of equations, say $(\ast)$, homogeneous separately in each set of variables $p$ and $\ph$, linear in  $\ph_{ij}^l$  and quadratic in  $p_{i_1\dotsc i_k}$. By Proposition~\ref{T:e}, the subset $S'$ defined by the equations~(\ref{E:d})  and $(\ast)$ is closed in $M$. Therefore $S' \cap A_{1\dotsc k}$ is closed in $A_{1\dotsc k}$. But $S \cap A_{1\dotsc k}$ equals $S' \cap A_{1\dotsc k}$ and hence is also closed.~$\bl$\\  

Now, consider the projections $\pi_1:S \to {\mathbb P}(\Lambda^k   {V})$ and $\pi_2: S \to {\mathbb P}(L)$ such that $ {\pi_1(p,\ph)=p}$, ${\pi_2(p,\ph)=\ph}$. These are regular maps. Clearly $\pi_1(S)=G(k,n)$, and $(\ph_1,\dotsc,\ph_t) $  belongs to $\pi_2(S)$  if and only if there is a vector subspace of dimension $k$ that is simultaneously isotropic for all $\ph_i$. Our lemma is equivalent to the assertion that $\pi_2(S)={\mathbb P}(L)$. By Proposition~\ref{T:b}, $ \pi_2(S)$ is closed in ${\Pm}(L)$, which is irreducible. Hence, by Proposition~\ref{T:B}, it is enough to show that $\dim \pi_2(S)=\dim{\mathbb P}(L)$ or just that $\dim \pi_2(S)\geq \dim{\mathbb P}(L)$.
 
\begin{lem} \label{T:q}
For any  $p\in G(k,n)$ the fibre $\pi_1^{-1}(p)$ is a projective variety of dimension
$$N=t\frac{n(n-1)-k(k-1)}{2}-1.$$ 
\end{lem}
\doc
Let $U$ be an arbitrary $k$-dimensional vector subspace of $V$. Choose a basis $\{ \itd e n \}$ of $V$ such that $\itd e k$ span $U$. Now $U$ is isotropic for an alternating bilinear form $\psi$ if and only if the matrix of $\psi$ has a zero $k \times k$ submatrix in the upper left-hand corner in this basis. The vector space consisting of  all such matrices has dimension  $\frac{n(n-1)-k(k-1)}{2}$.

For any point $p\in G (k,n)$ with $p=P(U)$ the fibre  $\pi_1^{-1}(p)$ consists of all  $t$-tuples of alternating bilinear forms, determined up to proportionality, such that $U$ is isotropic for all $\ph_i$. Consequently  
$\pi_1^{-1}(p)$ is a projective space of dimension $N$.~$\bl$

\begin{sled}
$S$ is irreducible variety of dimension
 \begin{equation} \label{E:J}
\dim S =  t\frac{n(n-1)-k(k-1)}{2}-1 + k(n-k).
\end{equation}
\end{sled}
\doc
Since $G(k,n)$ and a projective space are irreducible, we can apply Proposition~\ref{T:A} to   $S$ and $\pi_1$. It follows that $S$ is also irreducible. Hence, using Propositions~\ref{T:c}~{\sffamily(ii)} and~\ref{T:u}~{\sffamily(ii)}, we get
\begin{equation*}  
\dim S= \dim \pi_1^{-1}(p)+ \dim G(k,n)=  t\frac{n(n-1)-k(k-1)}{2}-1 + k(n-k).\ \bl
\end{equation*}

Choose a basis $\{\itd e n\}$ of $V$. Let $U$ be the vector subspace spanned by $\itd e k$. We denote by $L_U$  the subset of all $t$-tuples $(\itd \ph t) \in L$ such that $U$ is isotropic for all $\ph_i$, that is, 
\begin{equation} \label{E:n}
\ph^l_{ij}=0   \ \ \mbox{for } l\odo t \mbox{ and }  i,j\odo k.
\end{equation} 
It follows from Lemma~\ref{T:q} that
\begin{equation} \label{E:K}
\dim {\Pm}(L_U)=\dim \pi_1^{-1}(P(U))=N.
\end{equation}
 Put $S_U =\pi_2^{-1}({\Pm}(L_U))$, $\pi'_1=\pi_1|_{S_U}$, $\pi'_2=\pi_2|_{S_U}$.  Obviously, $S_U$ is a projective variety: $S_U$ is given by the system of equations that define $S$ and by the equations~(\ref{E:n}).
 
\begin{lem} \label{T:r}
Consider a point $p=P(U')\in G(k,n)$. Let  ${s= \dim U \cap U'}$. Then 
\begin{equation} \label{E:D}
\dim (\pi'_1)^{-1}(p)=  N-  \frac {kt(k-1)}2+\frac {st(s-1)} 2 . 
\end{equation}
\end{lem}
\doc
For any point $p \in G(k,n)$ with $p =P(U')$ the fibre $F=(\pi'_1)^{-1}(p)$ consists of the pairs $(p,\ph)$ such that $U$ and $U'$ are simultaneously isotropic for all of the forms of the $t$-tuple $\ph$.
Choose a basis  $\{\itd f n\}$ of $V$ such that $\itd f k$ span  $U$ and $f_{k-s+1},\dots,f_{2k-s}$ span $U'$. Now $U$ and $U'$ are isotropic for an alternating bilinear form $\psi$ if and only if the matrix of $\psi$ has two zero $k\times k$ submatrices on the main diagonal in this basis and these submatrices have common $s\times s  $  submatrix.
The vector space consisting of all such matrices has dimension $\frac{n(n-1)} 2 -2\frac {k(k-1)}2+\frac {s(s-1)} 2 $.  Therefore, since forms of $F$ are determined up to proportionality,
 
\begin{equation*}  
\dim F=t \left( \frac{n(n-1)} 2 -  {k(k-1)} +\frac {s(s-1)} 2  \right) -1 = N-  \frac {kt(k-1)}2+\frac {st(s-1)} 2 .\ \bl 
\end{equation*}

\begin{lem} \label{T:s}
If the numbers $k,\ t$ and $n$ satisfy the relations~(\ref{E:c}) and~(\ref{E:t2}), then  
\begin{equation} \label{E:o}
\dim S_U\leq N+ k(n-k)- \frac {kt(k-1)}2 .
\end{equation}
\end{lem}
\doc
Choose an irreducible component $S^0_U\subset S_U$ of maximal dimension, that is,\linebreak ${\dim S_U =\dim S^0_U}$. By Proposition~\ref{T:b}, $\pi'_1(S_U^0)$ is a projective variety.

Put $s_0= \max\{0, 2k-n\}$.  
Let $G_i=G_i(U,V)$ be closed subsets of $G(k,n)$ defined by~(\ref{E:Z}). Put $G_{k+1}=\varnothing $.
We have
$${G(k,n)=G_{s_0}\supset G_{s_0+1} \supset  \dots \supset G_k=\{P(U)\}\supset G_{k+1}=\varnothing}.$$
It follows from this that there exists a unique $s$ such that 
   $\pi'_1(S_U^0) \subset G_{ s}$ and   $\pi'_1(S_U^0) \not\subset G_{ s+1}.$ Then $\pi'_1(S_U^0) \cap({G_s \setminus G_{s+1}})\ne \varnothing$, and hence there exists $p\in \pi'_1(S_U^0) $ with $p=P(U')$ such that ${\dim (U\cap U')=s}$.
By Proposition~\ref{T:c}~{\sffamily (i)}, we obtain 
\begin{align*}
 \dim S^0_U &\leq\dim({\pi'_{1}}^{-1}(p)\cap S_U^0) +\dim \pi'_1(S_U^0) \leq\\
&\leq\dim{\pi'_{1}}^{-1}(p) +\dim G_s \overset {(\ref{E:L}), (\ref{E:D})}{=} \\
&=N- \frac {kt(k-1)}2+\frac {st(s-1)} 2 + (k-s)(n-k+s)=
 \\&
=s^2 \left( \frac t 2 -1 \right)+s\left(2k-n-\frac t 2  \right) -\frac{kt(k-1)}2 +k(n-k)  +N. 
\end{align*}

Thus,
$$
\dim S_U\leq N+\max_{s_0\leq s\leq k}  f(s), 
$$
where $f(s)= s^2 \left( \frac t 2 -1 \right)+s\left(2k-n-\frac t 2  \right) -\frac{kt(k-1)}2 +k(n-k) $.

If $t\geq 2$, then $\frac t 2 -1\geq 0$, and the function $f(s)$ attains its maximum either at  $s=s_0$ or at $s=k$. 

Note that the relations~(\ref{E:c}) and~(\ref{E:t2}) give
$$2k-n\leq 1 .$$
Indeed, suppose that $n< 2k-1$. Then, using~(\ref{E:c}), we get
$$ 4k-2 > t(k-1)+2k.$$ 
It is easy to see that this inequality holds true only if $t<2$. This contradicts~(\ref{E:t2}).

Therefore either $s_0=0 $ or $s_0=2k-n=1.$  Any way,
$$
f(s_0)= k(n-k)- \frac {kt(k-1)}2.
$$
We  also have
$$
f(k)=0.
$$

It follows from~(\ref{E:c}) that $f(s_0)\geq 0$. 
 Consequently $\max_{s_0\leq s\leq k}  f(s) =f(s_0)$, and
\begin{equation*}  
 \dim S_U\leq N+f(s_0). \ \bl
\end{equation*}

\begin{lem}\label{T:t}
If the numbers $k,\ t$ and $n$ satisfy the relations~(\ref{E:c}) and~(\ref{E:t2}), then there exists $\ph\in {\Pm}(L)$
such that\begin{equation}\label{E:p}
\dim \pi_2^{-1}(\ph)\leq k(n-k)- \frac {kt(k-1)}2 .
\end{equation}
\end{lem}
\doc
Let $S_U=\bigcup_j S_U^j$ be the decomposition of $S_U$ into irreducible components. By Proposition~\ref{T:b}, all the sets $\pi'_2(S_U^j)$ are closed in ${\Pm}(L_U)$. Hence we can assume without loss of generality that $\pi'_2(S^j_U)={\Pm}(L_U)$ for $j\leq j_0$ and for some $j_0$, and that $\pi'_2(S^j_U)$ is a proper closed subset of ${\Pm}(L_U)$
for $j>j_0$.  Since $\pi'_2(S_U)={\Pm}(L_U)$ and ${\Pm}(L_U)$ is irreducible,  we have $j_0\geq 1$.  It follows from Proposition~\ref{T:c}~{\sffamily (ii)} that for any  $j\leq j_0$  there exists a nonempty open subset $L_j\subset{\Pm}(L_U)$  such that for any point  $\ph\in L_j$
\begin{align*}
\dim ((\pi'_2)^{-1}(\ph) \cap S_U^j)& = \dim S_U^j- \dim \pi'_2(S_U^j) \leq \dim S_U  -  \dim {\Pm}( L_U) \overset{(\ref{E:K}),(\ref{E:o})}{\leq}\\ 
& \leq k(n-k)- \frac {kt(k-1)}2 .
\end{align*}
The intersection the finite number of nonempty open subsets of an irreducible set is always a nonempty set. Therefore $L_0=\bigcap_{j\leq j_0} L_j \setminus (\bigcup_{j>j_0} \pi'_2(S_U^j))$ is not empty, and for any $\ph \in L_0$ we have
\begin{equation*}  \dim \pi_2^{-1}(\ph)=\dim (\pi'_2)^{-1}(\ph)=\dim \bigcup_j ((\pi'_2)^{-1}(\ph) \cap S_U^j) \leq k(n-k)- \frac {kt(k-1)}2 .\ \bl \end{equation*}

We recall that, to prove the Main Lemma, it is enough to show that 
  $\dim \pi_2(S) \geq \dim {\Pm}(L)$. Chose a point $\ph\in {\Pm}(L)$ in accordance with Lemma~\ref{T:t}. Then, using Proposition~\ref{T:c}~{\sffamily (i)}, we get 
\begin{align*}
\dim \pi_2(S)&\geq \dim S-\dim \pi_2^{-1}(\ph) \overset{(\ref{E:J}),(\ref{E:p})}{\geq }\\ &\geq t\frac{n(n-1)-k(k-1)}{2}-1 + k(n-k)-  \left(k(n-k)- \frac {kt(k-1)}2 \right)= \\ &=t \frac {n(n-1)} 2 -1 = \dim {\Pm}(L).
\end{align*}
This concludes the proof of the lemma.

\section {Proof of the Theorem} \label{s:5}
\subsection{Reduction of the Theorem to the case of nilpotent Lie algebras of class 2.} \label{s:6}
\begin{utv} \label{T:f}
If $K$ is the complex or real   field, then the following equality holds:
$$g^n_{K}(s)=l^n_{K}(s).$$
\end{utv}
\doc
Let $G$ be a nilpotent Lie group of class 2 such that the dimensions of all abelian Lie subgroups of $G$ are not greater than $s$, and let  ${\mathfrak g}$ be the Lie algebra associated with $G$, which is also nilpotent of class 2. 
Consider an abelian subalgebra ${\mathfrak h}\leq{\mathfrak g}$. We denote by ${\mathfrak h}^M$ the minimal subalgebra such that ${\mathfrak h}\leq{\mathfrak h}^M$  and  there exists  a connected Lie subgroup $H\leq G$   with Lie algebra ${\mathfrak h}^M$.
By Theorem~1.4.3 of (Vinberg, Onishchik,~1995), the commutator subalgebras of ${\mathfrak h}$  and ${\mathfrak h}^M$  are equal. Therefore ${\mathfrak h}^M$ is commutative and hence so is $H$. By assumption,  ${\dim {\mathfrak h}\leq \dim {\mathfrak h}^M=\dim H\leq  s}$. We see that  ${\mathfrak g} $ has no commutative subalgebras of  dimension  greater than $s$. Consequently, by the definition of  $l^n_{K}$, 
$\dim G=\dim {\mathfrak g} \leq l^n_{K}(s)$. It follows that
$$g^n_{K}(s)\leq l^n_{K}(s).$$

On the other hand, given a complex (or real) nilpotent Lie algebra ${\mathfrak g}$ of class 2, there exists a connected complex (respectively real) nilpotent Lie group $G$ of class 2 with this Lie algebra (see Theorem~6.2 of (Vinberg, Onishchik,~1995)). 
A commutative Lie group  has a commutative Lie algebra. Therefore if
all commutative subalgebras of ${\mathfrak g}$  have dimensions at most $s$, then $G$ has no commutative subgroups of dimension greater than $s$.  
  This gives
$$g^n_{K}(s)\geq l^n_{K}(s).~\bl$$

Consider a nilpotent Lie algebra ${\mathfrak g}$ of class 2 over a field $K$. Let ${\mathfrak z}$ be the centre of  ${\mathfrak g}$,  $\{\itd{z}{t}\}$  a basis of  ${\mathfrak z}$,  and $V$ a complementary subspace to  ${\mathfrak z}$ of dimension $n$. Then the product of two elements  $x=\Bar{x}+\Bar{\Bar{x}}$ and $y=\Bar{y}+\Bar{\Bar{y}}$ of ${\mathfrak g}$     with $\Bar{x},\Bar{y}\in V $ and $ \Bar{\Bar{x}},\Bar{\Bar{y}}\in{\mathfrak z}$ has the form
\begin{equation} \label{E:g}
 [x,y]=\ph_1(\Bar{x},\Bar{y})z_1+ \dots+ \ph_t(\Bar{x},\Bar{y})z_t    
\end{equation}
for some $t$-tuple of alternating bilinear forms ${\Phi=\Phi({\mathfrak g})=(\itd{\ph}{t})}$ on $V$.

On the other hand, given vector spaces   ${\mathfrak z}$ and $V$ over $K$, a basis  $\{\itd{z}t\}$ of ${\mathfrak z}$, and a $t$-tuple ${\Phi=(\itd \ph t)}$ of alternating bilinear forms on $V$, one can define the product of two elements of $
{\mathfrak g}={\mathfrak g}(\Phi)= {\mathfrak z} \oplus V$ by~(\ref{E:g}). Obviously, ${\mathfrak g}(\Phi)$ is a nilpotent Lie algebra of class 2 with central subalgebra~${\mathfrak z}$.

 In this notation, $x$ and $y$ commute if and only if
\begin{equation*}  
 \ph_i(\Bar{x},\Bar{y})=0\ \mbox{ for } i=1,\dots, t. 
\end{equation*}
So that a subalgebra  $\mathfrak h \leq g$ is commutative if and only if a vector subspace   ${\mathfrak h}/ ({\mathfrak h}\cap {\mathfrak z}) \subset V$ is simultaneously isotropic for all of the forms  $\itd{\ph}{t}$. 

  Nilpotent associative algebras of class 2 have similar structure. And so we have the following proposition.
 
\begin{utv} \label{T:g}
For any field $K$
$$a^n_{K}(s)=l^n_{K}(s).$$
\end{utv}
\doc 
Let $A$ be a nilpotent associative algebra of class 2, and let $Z$ be its annihilator, that is, $Z=\{ x\in A \mid Ax=xA=0\}$. Further, let $\{\itd z t\}$ 
be a basis of  $Z$, and $V$ a complementary subspace to $Z$ of dimension $n$.  Then the product of two elements  $x=\Bar{x}+\Bar{\Bar{x}}$ and $y=\Bar{y}+\Bar{\Bar{y}}$ of $A$ with $\Bar{x},\Bar{y}\in V $ and $ \Bar{\Bar{x}},\Bar{\Bar{y}}\in Z$ has the form
\begin{equation} \label{E:h}
 x \centerdot y =\psi_1(\Bar{x},\Bar{y})z_1+ \dots+ \psi_t(\Bar{x},\Bar{y})z_t    
\end{equation}
for some bilinear forms $\itd \psi t$ on $V$.

And, as above, to any $t$-tuple of bilinear forms  $\Psi=(\itd \psi t)$ on a vector space $V$ one can assign a nilpotent associative algebra $A =A(\Psi)$ of class 2, where $A(\Psi)$ is a direct sum of $V$ and $Z$, and 
the product is given by~(\ref{E:h}).

Here elements $x$ and $y$ commute if and only if
\begin{equation*} 
 \psi_i(\Bar{x},\Bar{y})=\psi_i(\Bar{y},\Bar{x})\ \mbox{ for } i=1,\dots, t. 
\end{equation*}

Let $\itd C t$ be the matrices of $\itd \psi t$ with respect to some basis of  $V$.
 Then $x$ and $y$ commute if and only if 
$$
 \Bar x^\top C_i \Bar y= \Bar y^\top C_i \Bar x= \Bar x^\top C_i^\top \Bar y \ \mbox{ for } i=1,\dots, t,
$$
that is, if and only if
\begin{equation} \label{B}  
\Bar x^\top( C_i - C_i^\top )\Bar y=0\ \mbox{ for } i=1, \dots, t.
\end{equation}
Let $\itd {\psi^-} t$ be alternating bilinear forms given by matrices ${( C_1 - C_1^\top )}, \dots, {( C_t - C_t^\top )}$. It follows from~(\ref{B}) that  a subalgebra  $B \leq A$ is commutative if and only if a vector subspace  ${B / (B \cap Z) \subset V} $ is simultaneously isotropic for all of   $  \psi^-_i $.

Thus for any nilpotent associative algebra $A$ of class 2  the nilpotent Lie algebra ${\mathfrak g}(\Phi)$ of class 2 with ${\Phi=( \psi^-_1(A), \dots, \psi^-_t(A))}$ has the same dimension, and the dimensions of maximal commutative subalgebras of $A$ and ${\mathfrak g}(\Phi)$ are equal.

Conversely, for any nilpotent Lie algebra ${\mathfrak g}$  of class 2 there exists a nilpotent associative algebra $A(\Psi)$ of class 2 and of the same dimension such that the dimensions of maximal commutative subalgebras of  ${\mathfrak g}$ and $A(\Psi)$ are equal. Indeed, it is enough to choose $\Psi$ such that $ \psi_i^-=\ph_i({\mathfrak g})$.  
For example, put   $C^k_{ij}=\ph^k_{ij}({\mathfrak g})$ if $i<j$  and $C^k_{ij}=0$ otherwise.

Whence it follows, by Definition, that $a_K^n\equiv l_K^n.$~$\bl$

It follows from Propositions~\ref{T:f} and~\ref{T:g} that, to conclude the proof of the Theorem, it is enough to evaluate the function $l_K^n$.

\subsection{Proof of the Theorem for the function $l_K^n(n)$.} \label{s:7}

\begin{lem} \label{T:i}
Let $K$ be an algebraically closed field. Then the function $l_K^n$  satisfies the inequality 
$$
  l_K^n(s)\leq \max  \left( \frac{s^2+4} 8 +s, 2s-1  \right).
$$
\end{lem}
\doc
Consider a nilpotent Lie algebra  ${\mathfrak g}$ of class 2 over $K$. We   use the notations of the previous section.
 
Suppose that  $t\geq 2$.  Put $k=\left[ \frac{2n+t}{t+2} \right]$. Then  the relation~(\ref{E:c}) holds. It follows from the Main Lemma that there exists a $k$-dimensional subspace of $V$ that is simultaneously isotropic for all of the forms $\itd{\ph}t$. Consequently ${\mathfrak g}$ contains an abelian subalgebra of dimension   $s=k+t$. We have
$$
k\geq \frac {2n+t-(t+1)}{t+2}=\frac{2n-1}{t+2}.
$$
Thus we get
\begin{equation*}  
\dim{\mathfrak g}=n+t \leq \frac{k(t+2) +2t+1}{2}=\frac{t(s-t)+1} 2 +s \leq \frac{s^2+4} 8 +s.
\end{equation*}

If  $t=1$, then the form  $ \ph_1$ is nondegenerate, and there exists a basis $\{\itd e n\}$ for $V$ such that  $\ph_1$ is represented in this basis as follows:
 
\begin{equation} \label{E:k}
\ph_1=
\begin{pmatrix}
            0 & 1 & \cdots &  0 & 0 \\
          - 1 & 0 & \cdots & 0 & 0\\
          \hdotsfor{5} \\
           0 & 0 & \cdots &   0 & 1\\   
           0 & 0 & \cdots &  -1 & 0  
     \end{pmatrix} .
\end{equation}
Hence n is even, and  ${\mathfrak g}$ is isomorphic to the Heisenberg algebra $ {\mathfrak h}_{\frac n 2}$ of dimension $n+1$. Obviously,
the subalgebra spanned by $z_1$ and by all vectors $e_i$ with even $i$ is abelian. It has dimension  $s=\frac n 2 +1$. We have
$$
\dim {\mathfrak g} = 2s-1. 
$$

This concludes the proof.~$ \bl$

\begin{lem} \label{T:l}
For any field $K$ the following inequality holds:
$$
l^n_K(s)\geq 2s-1.
$$  
\end{lem}
\doc
We recall that if $\mathfrak g$ is the Heisenberg algebra  $ {\mathfrak h}_m$  of dimension $2m+1$, then $t=1$ and  $\ph_1$  is given by~(\ref{E:k}).

Consider elements $\itd x l, \ y\in V$  with coordinates   $ (\itd{x^1}{2m}), \dots,$ $ (\itd{x^l}{2m}),$ $(\itd y {2m})$ respectively. The element
 $y$ commutes with any  $x_i$ with $i\odo l$ if and only if
$$ \ph_1(x_i,y)=0\ \mbox{ for } i\odo l.$$
This is equivalent to the following  system of   equations:
$$\begin{matrix}
x_2^1y_1-x_1^1y_2+ \dots + x_{2m}^1y_{2m-1}-x_{2m-1}^1y_{2m}=0,\\
\hdotsfor{1}\\
x_2^ly_1-x_1^ly_2+ \dots + x_{2m}^ly_{2m-1}-x_{2m-1}^ly_{2m}=0.
\end{matrix}
$$

We notice that this system of   equations (in $y$) is linearly independent if and only if elements $\itd x l$ are linearly independent. 

Therefore, if  $\itd x m\in V$ are linearly independent elements such that $[  {x_i}, {x_j}]=0$ for $i,j\odo m$, then the set of all   $y\in V$  that commute with any $x_i$ is the vector subspace spanned by $\itd x m$.
Consequently  ${\mathfrak h}_m$   has no abelian subalgebras of dimension greater than  $s=m+1$. Note that
$\dim {\mathfrak h}_m=2s-1$.
The lemma follows from this.~$\bl$

\begin{lem} \label{T:m}
If $K$ is an infinite field, then
 $$
  l_K^n(s)\geq   \frac{s^2-1} 8 +s .
$$
\end{lem}
\doc
  Theorem~4$'$ of (Milentyeva,~2006) asserts that for any field $K$ the following inequality holds:
$$l_K^n(s)\geq \frac{s^2+4s-5} 8 .$$
If $K$ is infinite, the proof of this theorem can be modified to get a better estimate. Indeed, it was shown in Lemma~11 of (Milentyeva,~2006) that if the positive integers  $k,\ t$ and $n$ satisfy the inequality 
\begin{equation} \label{E:cc}
2n< t(k-1),
\end{equation}
and $V$ is an $n$-dimensional vector space over $K$, then there exists a $t$-tuple   $\Phi=(\itd \ph t)$ of alternating bilinear forms on $V$ such that no $k$-dimensional subspace is simultaneously isotropic for all of the forms  $\itd \ph t$.
The proof of that lemma consists of two cases: the case where $K$ is a finite field, and the case where $K$ is infinite. 
In the infinite case,  the Main Lemma of (Milenteva,~2004) is used. In that lemma the following weaker inequality is sufficient:
\begin{equation} \label{E:bb}
2n< t(k-1)+2k.
\end{equation}
So, in the infinite case, the condition~(\ref{E:cc}) of Lemma~11 of (Milentyeva,~2006) can be replaced by~(\ref{E:bb}).
 (This   can easily be  checked for algebraically closed field $K$ using the proof of our Main Lemma. Indeed,
${\dim \pi_2(S) \leq \dim S < \dim {\Pm}(L)}$
  whenever the relation~(\ref{E:bb}) holds.)
 The arguments of \S\ref{s:6} show  that   the corresponding Lie algebra ${\mathfrak g}(\Phi)$ has no abelian subalgebras of dimension greater than $k+t-1$.
 
Let $s$ be even. Put  $ t=\frac s 2,\ k=\frac s 2 +1,\      n=\left[ \frac{s^2+4s+7}{8}\right] $. Then ${2n< t(k-1) +2k}$ and, by the previous arguments, there exists a Lie algebra ${\mathfrak g}(\Phi)$ all of whose abelian subalgebras have dimension at most $s$ and such that the dimension of ${\mathfrak g}(\Phi)$ is equal to 
\[n+t=\left[\frac{s^2+4s+7}{8}\right]+\frac{s}{2} \geq\frac{s^2}{8}+s .\]

Similarly, if $s$ is odd, putting  $t=k=\frac{s+1}{2}$,  $n=\left[ \frac{s^2+4s+2}{8}\right] $, we obtain that the dimension of ${\mathfrak g}(\Phi)$
equals 
\[  n+t=\left[\frac{s^2+4s+2}{8}\right]+\frac{s+1}{2} \geq\frac{s^2-1}{8}+s .\]

Hence   $ l_K^n(s) \geq  \frac{s^2-1}{8}+s $ and the Lemma is proved.~$\bl$\\

Note that  $l_K^n$ is an integer function, and that the difference between $\frac{s^2+4} 8+s$ and $\frac{s^2-1} 8+s$ is less than 1. Combining Lemmas~\ref{T:i},~\ref{T:l}, and~\ref{T:m}, we get that if $K$ is   algebraically closed, then 
 $$
  l_K^n(s)= \max  \left( \left[\frac{s^2+4} 8\right] +s, 2s-1  \right).
$$

This implies the relation~(\ref{E:m}), and the Theorem is proved.\\

\textbf{Acknowledgment. } The author would like to express her gratitude to Professor A.Yu.~{Olshanskii}  for  suggesting the problem,  his constant attention to this work, and his valuable advice.

{\scriptsize Higher Algebra, Department of Mechanics nd Mathematics, M.V.~Lomonosov Moscow State University, Leninskie Gory, GSP-2, Moscow, Russsia; 119992. }

{\scriptsize  E-mail: mariamil@yandex.ru}

\end{document}